\documentclass[letterpaper, 10 pt, conference]{ieeeconf}  % Comment this line out
\IEEEoverridecommandlockouts                              % This command is only
\usepackage[utf8]{inputenc}
\usepackage[T1]{fontenc}    % use 8-bit T1 fonts

\usepackage{hyperref}       % hyperlinks
\usepackage{amssymb,amsmath,amsfonts}
\usepackage[all,arc,cmtip]{xy}
\usepackage{tikz}
\usepackage{enumerate}
\usepackage{mathrsfs}
\usepackage{algorithm}
\usepackage{graphicx}
\usepackage{booktabs}
\usepackage{longtable}
\usepackage{wrapfig}
\usepackage{pdflscape}
\usepackage{nicefrac}       % compact symbols for 1/2, etc.
\usepackage{microtype}      % microtypography

\usepackage{dcolumn}

\newcommand{\tr}{\operatorname{tr}}

\usepackage{microtype}
\usepackage{graphicx}
\usepackage{subfigure}
\usepackage{booktabs} % for professional tables

\newtheorem{proposition}{Proposition}

\author{}
\title{\LARGE \bf A Two-Step Pre-Processing for Semidefinite Programming}

\author{Vyacheslav Kungurtsev and Jakub Marecek
\thanks{Equal contribution. Both authors at the Department of Computer Science,
  Czech Technical University,
  Prague, the Czech Republic.}
\thanks{\small \tt vyacheslav.kungurtsev@fel.cvut.cz}
\thanks{\small \tt jakub.marecek@fel.cvut.cz}
}

\begin{document}

\maketitle
\thispagestyle{empty}
\pagestyle{empty}


%\maketitle

\begin{abstract}
In semidefinite programming (SDP), a number of pre-processing techniques have been developed, including procedures based on chordal decomposition, 
which exploit sparsity in the semidefinite program in order to reduce the dimension of individual constraints, and procedures based on facial reduction, which reduces the dimension of the problem by removing redundant rows and columns.  
So far, these have been studied in isolation.
We show that these techniques are, in fact, complementary. 
In computational experiments, we show that a two-step pre-processing followed by a standard interior-point method outperforms the interior point method, with or without either of the pre-processing techniques, by a considerable margin.
%For example, SeDuMi benefiting from the two-step pre-processing solves 50\% of a subset of SDPLib, where it is applicable, within 0.19 seconds, whereas without the pre-processing solving 50\% of the same subset takes , in both cases to the same tolerances on the same hardware.
\end{abstract}

\section{Introduction}

There has been much recent interest in semidefinite programming (SDP),
based on the realisation that it provides very tight relaxations for
 non-convex problems in a variety of domains, including 
statistics and core machine learning  \cite{NIPS20042628,dahl2008covariance,NIPS20176645}, 
computer vision \cite{schweighofer2008globally},
automatic control \cite{boyd1994linear,7979611},
and robotics \cite{6907313}.
Often, these instances can be seen as relaxations of  certain non-convex polynomial optimisation problems \cite{schweighofer2008globally,Lavaei2012,Waki2012,NieDemmel,NIPS20176645}.

Instances of semidefinite programming obtained as relaxations of polynomial optimisation problems, among others, are 
sparse, structured, and strictly feasible (cf. Theorem 3.2 in \cite{NieDemmel}).
One can exploit the structure and sparsity using the so-called chordal decomposition \cite{Grone,blair1993introduction,Fukuda2001,Kim2012,Bergman2016,zheng2020chordal}, which we introduce formally in the next section.
For example, one can use chordal-decomposition in pre-processing \cite{Fukuda2001}. 
Many examples of the use of chordal-decomposition pre-processing abound in control \cite{7979611},
power systems \cite{Marecek2016,6692873,6756976},
and statistics \cite{dahl2008covariance},
but it turns out that such pre-processing can lead to non-trivial numerical issues, as has been recently explained by \cite{Raghunathan2016}.

Other types of pre-processing can, in fact, address
numerical issues. Notably,  (partial) facial reduction \cite{borwein1981facial,Pataki2015,Pataki2017} addresses issues associated with a lack of a strictly feasible point, as defined in the next section.
Nevertheless, facial-reduction pre-processing has, so far, attracted much less interest, and has not been used in conjunction with the chordal decomposition, yet. 

Here, we introduce a two-step pre-processing for SDP, combining techniques from chordal decomposition and facial reduction.
We present extensive numerical results comparing the performance of a standard SDP interior-point method on its own, coupled with chordal-decomposition pre-processing, and coupled with the two-step pre-processing
on a variety of large-scale 
structured instances of SDP, including those arising from the venerable maximum bisection (MAXCUT) problem,
a benchmark in binary quadratic programming (BiqMac),
a benchmark in power-systems engineering (IEEE Test Cases),
and a benchmark in general-purpose semidefinite programming (SDPLib). 

Our key observation is that the two-step pre-processing appears 
%to perform no worse than either pre-processing routine on its own, or the interior-point method without any  pre-processing, and,
does in does improve the performance dramatically, when there is sparsity or structure.
Even on SDPLib, which is not know to exhibit any particular structure, considerable improvements are possible.
For example, following the two-step pre-processing, SeDuMi, a commonly used open-source interior-point method, solves 50\% of a subset of SDPLib, where it is applicable, more than 11 times faster than SeDuMi without pre-processing, in both cases to the same tolerances and on the same hardware and inclusive of the time it takes to perform the pre-processing.

\iffalse
\subsection{Summary of the Paper}
The rest of the paper proceeds as follows. We provide some background on semidefinite programming and the two pre-processing tools in Section~\ref{sec:back}. In Section~\ref{sec:alg} we present the two step pre-processing Algorithm. We present our numerical results in Seciton~\ref{sec:num} and summarize our findings and discuss potential future work in Section~\ref{sec:conc}.
\fi

\section{Background}\label{sec:back}
\subsection{Semidefinite programming}
\label{sec:SDP}

Let us recall some standard definitions. Consider an optimisation problem over the set $\mathbb{S}^n$ of symmetric $n \times n$ matrices:
\begin{align}\tag{SDP}\label{P}
\begin{array}{rl}
\min_{X\in\mathbb{S}^n} & C \bullet X \\
\text{subject to } & A_i \bullet X = b_i, \quad i = 1, 2, \ldots, m \\
& X\succeq 0
\end{array}
\end{align}
where 
$C, A_i \in \mathbb{S}^n$ are compatible matrices, $b_i \in \mathbb{R}^n$ are compatible vectors, 
$\bullet$ denotes the inner product on $\mathbb{S}^n$, i.e., $A \bullet B := \tr(A^{T}B)=\sum _{{i=1,j=1}}^{n}A_{{ij}}B_{{ij}},$ and $X\succeq 0$ denotes the constraint on matrix $X$ to be positive semi-define, i.e., $v^T Xv \ge 0$ for all $v \in \mathbb{R}^n$. This problem is known as semidefinite programming.
For simplicity, we assume that it is feasible.

As usual, for any convex set $C \subseteq \mathbb{S}^n$, a point $X \in \mathbb{S}^n$ lies in the relative
interior of $C$ if and only if
$\forall y \in C, \exists z \in C, \exists \alpha \in R, 0 < \alpha < 1: x = \alpha y + (1-\alpha) z$. 
The relative interior of the cone of positive semidefinite matrices are the positive definite matrices, i.e., $X$ where $v^T Xv > 0$ for all $v \in \mathbb{R}^n \setminus \{ 0 \}$.

\subsection{Chordal decomposition}

Chordal decomposition is an established pre-processing technique in  semidefinite programming, with a history of research going back to 1984 \cite{Grone},
with a considerable revival \cite{Fukuda2001,vandenberghe2015chordal,7963462,Raghunathan2016,zheng2016cdcs,zheng2019chordal,zheng2020chordal} in the past two decades. % andersen2010linear
It is also known as Matrix Completion Pre-processing 
or the d-space and r-space Conversion Method
\cite{Fukuda2001,Kim2012,Bergman2016}.
There are important applications in
statistics \cite{dahl2008covariance},
machine learning \cite{wainwright2004treewidth,wainwright2008graphical,NIPS20176645},
power systems \cite{Marecek2016,6692873,6756976}, and
automatic control \cite{7979611}.
% For notable examples of the use of the decomposition, we refer to \cite{7979611,6756976,Marecek2016,6692873,dahl2008covariance}. 
%Within machine learning, related techniques have been 
%used in graphical models 
%\cite{wainwright2008graphical,wainwright2004treewidth,NIPS20176645}.
For an excellent survey, see \cite{vandenberghe2015chordal}.

Let us define chordal decomposition of \eqref{P} with structured  $A_1,\ldots,A_m \in S^n$ formally.
First, let us define the (so-called correlative) sparsity pattern as a simple undirected graph $G=(N, E)$, where $N:=\{1, \ldots, n\}$ and

\begin{equation*} E : =\{(j, k):j\neq k, [A_{i}]_{jk}\neq 0\,\,\, \text{ for some } i=1, \ldots, m\}. \end{equation*}

Given the sparsity pattern $G(N,E)$ of an SDP, chordal decomposition computes a chordal extension $F$ with $E \subseteq F$, a set of maximal cliques ${C_1,\ldots,C_l}$ of the graph $G(N,F)$, and a clique tree
$\mathcal{T}(\mathcal{N},\mathcal{E})$. Using a mapping $\sigma_s:\mathbb{N}\to\{1,...,|C_s|\}$ (where by $|A|$ we mean the cardinality of set $A$) from the original indices to an ordering of the clique $C_s$, we can define:

\begin{align} 
[A_{s, p}]_{\sigma_{s}(i)\sigma_{s}(j)} & = \begin{cases} [A_{p}]_{ij} & \text{if}\ s=\min\{t\vert (i, j)\in \text{C}_{t}\}\\ 0 & \text{otherwise} \end{cases} \label{eq:reorder1}\\ 
E_{s, ij} & = \frac{1}{2}\left(e_{\sigma_{s}(i)}e_{\sigma_{s}(j)}^{T}+e_{\sigma_{s}(j)}e_{\sigma_{s}(i)}^{T}\right)\forall i, j\in \text{C}_{s} \label{eq:reorder2}\\ 
(s, t) & \in \mathcal{T}\Leftrightarrow(\text{C}_{s}, \text{C}_{t})\in \mathcal{E} \label{eq:reorder3}\\ 
\text{Q} & = \text{C}_{s}\cap \text{C}_{t} \label{eq:reorder4}
\end{align}
with the notation that $e_{\sigma_s(i)}\in \mathbb{R}^{|C_s|}$.

The SDP is then reformulated using a positive semi-definiteness constraint for each maximal cliques and equality constraints for any vertices in more than one maximal clique:

\begin{equation} \begin{split} \min_{X_{s}\in \mathbb{S}\vert \text{C}_{s}\vert} & \sum_{s=1}^{\ell}A_{s,0}\bullet X_{s}  \label{MC}\\ \text{subject to} & \sum_{s=1}^{\ell}A_{s, p} \bullet X^{s}=b_{p}\qquad \forall p=1, \ldots, m\\ & E_{s, ij}\bullet X_{s}=E_{t, ij}\bullet X_{t}\quad \forall i\leq j, i, j\in \text{Q},\\ &\qquad\qquad\qquad\qquad\qquad (s, t)\in \mathcal{E}\\ & X_{s}\succeq 0\qquad\qquad\qquad\quad \forall s=1, \ldots, \ell. \end{split}  \end{equation}

\subsection{Facial reduction}

There is also a long history of work on facial reduction
\cite{kriswolk:09,Cheung:2013,7040427,ChDrWo:14,XinghangYe:16,permenter2017reduction,Permenter2017,hu2019facial}, 
another type of pre-processing for SDPs.
%including some of the most elegant papers in optimisation \cite{kriswolk:09}. 
\iffalse
The conic expansion method \cite{luo1997duality}
is in some sense the dual of facial reduction.
\fi
Facial reduction has been used to pre-process \cite{MR3108430,DrusLiWolk:14} degenerate semidefinite programs.
Its applications to machine learning have been limited \cite{Babakmanifold:12,HuangWolkXYe:16,RWW:17,rpca} so far.
For a quick overview of facial reduction and its relationship to degeneracy, we refer to \cite{DrusWolk:16}.

To define facial reduction formally, let us present some facts from convex geometry following \cite{DrusWolk:16}. Let us consider a convex cone $\mathcal{K}$ such as $\mathbb{S}^n$. 
A convex cone $F \subset \mathcal{K}$ is called a face of $\mathcal{K}$, when $x, y \in \mathcal{K}, x + y \in F \implies x, y \in F$.
A face of $\mathcal{K}$ is called proper, if it is neither empty, nor $\mathcal{K}$ itself.
Clearly, the intersection of an arbitrary collection of faces of $\mathcal{K}$ is itself a (possibly lower-dimensional) face of $\mathcal{K}$.
Non-trivially, the relative interiors of all faces of $\mathcal{K}$ form a partition of $\mathcal{K}$, i.e., every point in $\mathcal{K}$ lies in the relative interior of precisely one face and any proper face of $\mathcal{K}$ is disjoint from the relative interior of $\mathcal{K}$.

Next, let us consider the
dual of $\mathcal{K}$, and denote it by $\mathcal{K}^*$.
We use $\mathcal{K}^\perp$ for the orthogonal complement of the affine hull of $\mathcal{K}$.
Any set of the form $F = v^\perp \cap \mathcal{K}$, for some $v \in \mathcal{K}^*$, is called an \emph{exposed face} of $\mathcal{K}$ with an exposing vector $v$. 
It is well-known (Proposition 2.2.1 in
\cite{DrusWolk:16})
that if faces $F_1, F_2 \subset \mathcal{K}$ are exposed by vectors $v_1, v_2 \in \mathcal{K}^*$, then the intersection $F_1 \cap F_2$ is exposed by $v_1 + v_2$.
Finally, a convex cone is exposed if all its faces are exposed.

In the case of $\mathbb{S}^n$, and its dual $\mathbb{S}^n$ (self-duality), there is a correspondence between $r$-dimensional linear subspaces $\mathcal{R}$ of $\mathbb{R}^n$ and faces of $\mathbb{S}^n$, wherein 
$F_{\mathcal{R}} := \{ X \in \mathbb{S}^n \textrm{ such that range } X \subseteq \mathcal{R} \}$ is a face of $\mathbb{S}^n$.
Consequently, for any matrix $V \in \mathbb{R}^{n \times r}$ with range $V$ equal to $\mathcal{R}$, we have $F_{\mathcal{R}} = V \mathbb{S}^r V^T$, i.e., the face 
is isomorphic to an $r$-dimensional positive semi-definite cone $\mathbb{S}^r$.
Subsequently, the $F_{\mathcal{R}}$ is being exposed by some $UU^T$ for $U \in \mathbb{R}^{n \times (n - r)}$.

In facial reduction, one considers an instance of \eqref{P}, where the Slater condition fails, and iteratively constructs an equivalent instance, which has a Slater point. 
In each iteration, one aims to find $y$ such that:
\begin{align}\tag{TEST}\label{test}
\sum_{i=0}^m y_i A_i \succeq 0 \textrm{ and } b \bullet y = 0.
\end{align}
If no such vector $y$ exists, the Slater condition holds by the Theorem of the Alternative (Theorem 3.1.2 in \cite{DrusWolk:16}).
Otherwise, we know that the minimal face containing the feasible set is contained in $\mathcal{K'} := (\sum_{i=0}^m y_i A_i)^{\perp} \cap \mathcal{K}$, which yields an equivalent:
\begin{align}
\tag{FR} \label{FR}
\min_{X\in\mathcal{K'}} & C \bullet X \\
\text{subject to } & A_i \bullet X = b_i, \quad i = 1, 2, \ldots, m \notag \\
& X \in \mathcal{K'} \notag
\end{align}

%Let us now remark on the use of facial reduction in case of violation of  Mangasarian-Fromovitz constraint qualification (CQ) rather than Slater CQ.
%This case has been studied by Cheung et al.  \cite[Section 4.1.2]{MR3108430}. 
%In principle, the violation of  Mangasarian-Fromovitz constraint qualification may require an additional  explicit treatment of certain implicit equality constraints  \cite[Proposition 3.23]{MR3108430}, which does not render the process substantially different and still allows for the recovery of a solution with a Slater point.
%In the following, we will focus on the case of the loss of Slater CQ for simplicity.
This leaves only the small matter of an efficient implementation.

\section{An Algorithm}\label{sec:alg}

\begin{algorithm}[t!]
%\framebox[3.35in]{
\parbox{3.05in}{ 
{\sc TwoStep}($C$, $A_i$, $b_i$)
\begin{tabbing}{}
******\=******\=*******\=***\= \hspace{2.5in} \=  \kill  

$(G(N,E),\mathcal{T}(\mathcal{N},\mathcal{E}))=$ {\sc Sparsity}($C, A_1, A_2, \ldots, A_m$) \\

$(A,E,Q)$ = {\sc Reorder} by~\eqref{eq:reorder1},~\eqref{eq:reorder2},~\eqref{eq:reorder3},~\eqref{eq:reorder4} \\ 

%$C,  A_2, \ldots, A_m$ = {\sc reorder}($C, A_1, A_2, \ldots, A_m$) \\ % fill-reducing ordering 
%$R$ = {\sc embed}($A_1, A_2, \ldots, A_m$) \\ % chordal embedding by elimination 
%$Q$ = {\sc cliques}($R$) \\ % chordal embedding by elimination 
%Gt = R+R'; \\
%Gt(p,p) = Gt; \\ % recover embedded pattern 
Let $F_0$ be the feasible set of \eqref{MC} \\
 Let $y_0 = 0, \, i=1.$ \\
 {\bf repeat }  \\
 \> {\bf if } no $y$ satisfying \eqref{test} with $F_{i - 1}$ exists, {\bf break } \\
 \> {\bf else } choose \= $y_{i} \in L \cap F_{i-1}^*.$                 \\
  \>        \= Compute  \eqref{SVD} \\
 \>        \= Let $F_i = F_{i-1} \cap y_i^\perp$, defined by \eqref{FR2}  \\
 \> Let  \= $i = i+1.$ \\
 {\bf end}    \\
{\bf return} $F_{i - 1}$
\end{tabbing}
}%}
\caption{A ``two-step'' pre-processing for an instance $C$, $A_i$, $b_i$ of semidefinite programming \eqref{P}.
%\JM{SLAVA: please feel free to edit. It needs some more clarity.}
} 
\label{alg:twostep}
\end{algorithm}
\vspace{.35cm}

In Algorithm~\ref{alg:twostep}, we suggest that our pre-processing has two steps, at a high level: chordal embedding and facial reduction.

In more detail, there are several substeps to the first step.
As a first substep, we may wish to compute a fill-reducing ordering of the matrices. 
In Matlab, for instance, function {\tt amd} can be used to compute an approximate minimum degree permutation vector for a sparse matrix,
which could be obtained as the union of support sets of the matrices $C$ and $A_i, i = 1, 2, \ldots, m$.
We note that this substep can be omitted, if necessary, without compromising the results of our analysis below.
As a second substep, we compute the chordal embedding.
In Matlab, for instance, function {\tt symbfact}  returns the sparsity pattern of the Cholesky factor as its fifth output, which is a widely used embedding. 
Based on the chordal embedding, 
one can list all maximal cliques by breadth-first search. 
That is, {\tt cliques} forms a clique for each vertex together with its neighbors that follow in a perfect elimination ordering, and tests whether the set of cliques is maximal.
Finally, we construct the new semidefinite program \eqref{MC} based on the cliques $C$.

Subsequently, we run the iterative facial reduction,
where we interweave substeps of 
testing whether to continue with vector $y$ obtained in \eqref{test}
and reducing the instance by computing the singular-value decomposition (SVD) of the positive semidefinite $\sum_{i=0}^m y_i A_i$ to obtain:
\begin{align}
\label{SVD}
\tag{SVD}
\begin{bmatrix}
U & V
\end{bmatrix}
\begin{bmatrix}
D & 0 \\
0 & 0
\end{bmatrix}
\begin{bmatrix}
U & V
\end{bmatrix}^T,
\end{align}
where $[U V] \in \mathbb{R}^{n_{i-1} \times n_{i-1}}$ is an orthogonal matrix 
with $n_{i-1}$ being the dimension of the $F_{i-1}$ in  the previous substep,
and $D \in \mathbb{S}^r$ is a diagonal matrix. The simplified $F_{i}$ is then:
\begin{align}
\label{FR2}
\tag{SDP-FR}
\min_{X'\in\mathbb{S}^{n_{i-1} - r}} & V^T C V \bullet X' \\
\text{subject to } & A_i \bullet V^T X' V = b_i, \quad i = 1, 2, \ldots, m \notag \\
& X'\succeq 0, \notag
\end{align}
where $n_{i-1}$ is, again, the dimension in the previous substep. Notice that in each iteration $r>0$, and hence there can be at most $n$ iterations.

It is not immediately obvious that algorithm 
{\sc TwoStep}, Algorithm~\ref{alg:twostep}, is efficient. Indeed, {\sc reorder} and {\sc embed} solve an NP-Hard problem \cite{Yannakakis1981} and we solve a number of non-trivial optimisation problems \eqref{test} and \eqref{FR}.

First, notice that algorithm {\sc TwoStep} does not require the minimum fill-in reordering or embedding:

\begin{proposition}[Based on \cite{Natanzon2000}]
\label{thm:embed}
There is an implementation of {\sc embed}, which 
in a graph $R$ with maximum degree $d$
and minimum fill-in $k$ 
produces a solution within a factor of $O(d^{2.5} \log^4(kd))$ of the optimum in time $O(knm + \min(n^2 M(k)/k, nM(n))$, where $M(n)$ denotes the number of operations needed to multiply two $n \times n$ Boolean matrices.
% https://www.researchgate.net/profile/Ron_Shamir/publication/2644373_A_Polynomial_Approximation_Algorithm_for_the_Minimum_Fill-In_Problem/links/0fcfd50e41c3e448bd000000/A-Polynomial-Approximation-Algorithm-for-the-Minimum-Fill-In-Problem.pdf?origin=publication_list
\end{proposition}

This is a reasonably tight result, considering that even a constant-factor approximation is NP-Hard, cf. Theorem 21 in \cite{Bodlaender1995},
exact algorithms cannot run within 
time $2^{O(\sqrt{n} / \log^c n) }$
assuming the Exponential Time Hypothesis,
%\cite{Bliznets2016},
%https://arxiv.org/pdf/1508.05282.pdf
and a variety of lower bounds \cite{Cao2018}.
We note that there are also algorithms \cite{Fomin2013} producing fill-in $k$, which is optimal, 
albeit not running in polynomial time.

Second, notice that:

\begin{proposition}
\label{thm:fr}
For a polynomial-time {\sc reorder} and {\sc embed}, Algorithm {\sc TwoStep} runs in  time polynomial in $n$ and $m$ on the BSS machine.
\end{proposition}

\begin{proof}
The proof is trivial, once one notices 
that any facial-reduction algorithm solves at most $n$ instances of \eqref{FR} in at most $n$ iterations of the seemingly infinite loop, as explained in Section 4.2 in \cite{DrusWolk:16}.
\end{proof}

Notice, however, that this result reasons about the behaviour of BSS machine \cite{Blum1997}, rather than the more usual Turing machine, due to the complexity of numerical routines involved, SVD or otherwise.

\section{Numerical results}\label{sec:num}

%\JM{Throughout, we should start with a comparison with solvers other than SeDuMi, which is very impressive, and only then present the comparison with SeDuMi.}

%The first step could, in principle, be run using {\tt SparseCoLo} of \cite{Fukuda2001,kim2011exploiting}.
%The second step could, likewise, 
%be run using  {\tt frlib} of \cite{Permenter2017}.
%Each of these packages is rather complex, though.

In this section we compare the reliability and speed of interior-point method SeDuMi without any pre-processing, with the chordal-decomposition pre-processing (implemented using SparseCoLO of \cite{fujisawa2009user,kim2011exploiting}),
and both chordal-decomposition and facial-reduction pre-processing (implemented using SparseCoLO and {\tt frlib} of \cite{permenter2014partial}). We chose SeDuMi as compared to other solvers due to its overall robustness~\cite{mittelmann2003independent} and the reported complementary performance of this algorithm with FR~\cite{permenter2017solving}.
The tests were performed on a computing cluster, using 4 cores running at 3 GHz and memory allocated as necessary,
%(e.g., for the maxcut problems, up to 300G, but for all others 50G) 
running MATLAB 2018b on Debian.

Our conjecture was that especially for sparse structured SDP,
chordal-decomposition improves the speed of convergence for some problems, while  producing problems too poorly structured for interior-point solvers to solve quickly and reliably. In contrast, the additional facial-reduction step corrects this, and ultimately results in an improved performance overall, compared to both other settings.

%The SeDuMi and MOSEK tests were performed on an Intel Core i7-4702MQ CPU @ 2.20 GHz x 8 with 7.7G of memory, running MATLAB 2015b on Ubuntu.
%cluster

We report some of the results as performance profiles. These were introduced in~\cite{Dolan2002} as a way of visualizing the dual
performance measures of robustness (solving the largest proportion of problems), and efficiency (solving them quickly). The level of each curve at the
right-hand vertical boundary indicates how many problems were solved, and the relative location of each curve compared to the others in the profile intermediately
indicates the speed of convergence compared to the best solver
for each problem. Simply put, the further to the upper left corner a curve corresponding to an algorithm is, the better.

\iffalse
Although we are aware of the fact that there are many alternatives \cite{Gould2016}, we still consider performance profiles a standard.
\fi

\subsection{SDPLib}

Our main experiment considers the SDPLib test set~\cite{borchers1999sdplib}, which is a standard for benchmarking SDP software~\cite{mittelmann2003independent}, composed of
a variety of toy, academic, and real-world SDP problems. These are known to be sparse, 
but no particular structure is shared across the test set. 
We refer to \url{http://plato.asu.edu/ftp/sparse_sdp.html} for details of the instances
and the results obtained by a variety of both free and commercial solvers.
Note that out of the 92 instances within the test set, 
leading solvers can solve 56--88 instances, within a 40000-second (11-hour) time limit per instance considered by \cite{mittelmann2003independent}. 

\begin{figure}
\begin{center}
\includegraphics[scale=0.6]{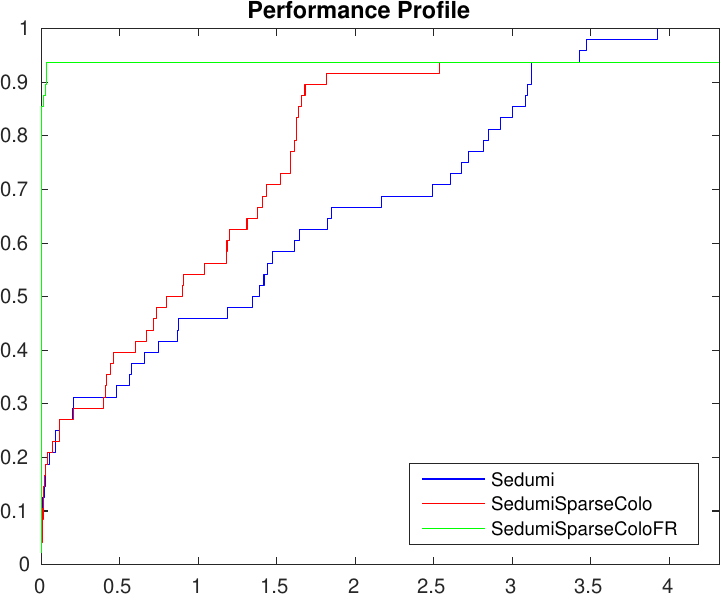}
\end{center}
\vskip 3mm
\caption{\label{fig-sdplibsedumi} Performance profiles (i.e., proportion of instances from a test set solved over time, relative to the best performer, cf. \cite{Dolan2002}) comparing SeDuMi with no (Sedumi), chordal-decomposition (SedumiSparseColo), and two-step processing (SedumiSparseColoFR) on 49 instances from SDPLib.}
\label{fig:sdplib}
\end{figure}

Rather surprisingly, 
we can demonstrate that one can improve the performance of SeDuMi, an open-source implementation of the interior-point method, 
which can work with instances without a Slater point.
%and altogether can create instances that are problematic for the interior point algorithms to solve.
%\subsubsection{SeDuMi Results}
Fig.~\ref{fig:sdplib} presents a performance profile \cite{Dolan2002} on 49 of the 92 instances, where running SparseCoLo did not report a failure. 
The corresponding numerical values are available from Tab.~\ref{table:times}: for instance, SeDuMi can solve 50\% of instances within 2.89 seconds without any pre-processing, 
or within 0.19 seconds with the two-step pre-processing.
In both cases, the run is to the same tolerances on the same hardware, inclusive of the time it takes to perform the pre-processing.
% We found that SparseCoLo resulted in a number of seg faults with MATLAB. 
Overall, we found that each step of pre-processing does 
improve the speed of convergence of the interior-point method in this test set. 
The reliability does seem to worsen, perhaps as a result of ill-conditioning of the linear systems employed in the interior-point method, after the chordal-decomposition pre-processing. 
Still, the proportion of the problems solved is very high among this test set,
and the strong increase in the speed due to the two-step pre-processing is an indication of the strength of the approach.

%\newcolumntype{d}[1]{D{.}{\cdot}{#1} }
\newcolumntype{d}{D{.}{.}{2} } %decimal column as before
%\newcolumntype{d}{D} %decimal column as before

\begin{table}[tb]
\centering
\caption{Run-time in seconds of SeDuMi with no, chordal-completion (w/ SC), and two-step (w/ SC+FR) processing, which has been required to solve a certain proportion of instances in SDPLib, up to the tolerances considered by Mittelmann \cite{mittelmann2003independent}. Dash (--) indicates the run did not finish.}
\vskip 3mm
\label{table:times}
\centering
\begin{tabular}{r d  d  d } 
\toprule
Proportion solved 
 & 
\multicolumn{1}{r}{SeDuMi} & 
\multicolumn{1}{r}{w/ SC} & 
\multicolumn{1}{r}{w/ SC+FR} \\  
%\cmidrule(r){1-5}
\midrule
$25\%$ & 0.48 & 0.23 & 0.084 \\
$50\%$ & 2.29 & 0.69 & 0.19 \\
$75\%$ & 5.07 & 2.79 & 2.8 \\
$94\%$ & 28 & 41 & 40 \\
$100\%$ & 43 & - & - \\
\bottomrule
\end{tabular}
\end{table}

\subsection{Polynomial optimisation}

\begin{figure}
\begin{center}
\begin{tabular}{c c c c} 
\includegraphics[scale=0.3]{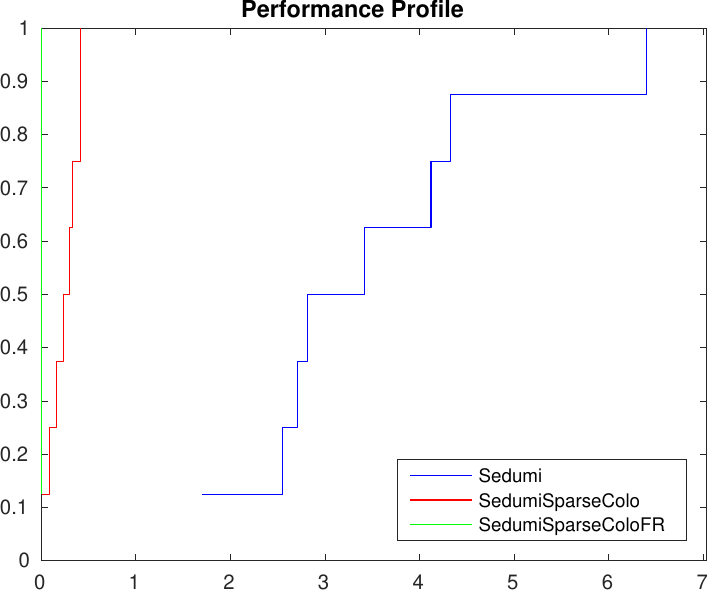}
&
\includegraphics[scale=0.3]{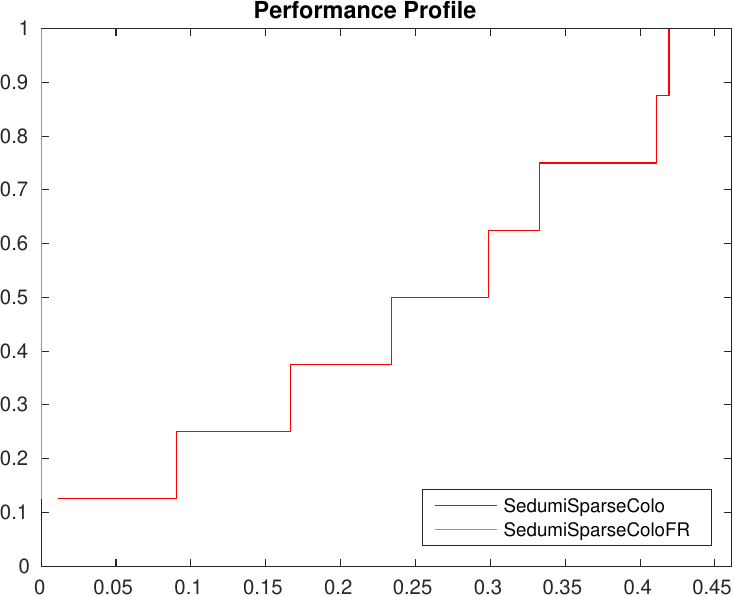} \\
Primal &
(Primal, zoomed-in)
\\
\includegraphics[scale=0.3]{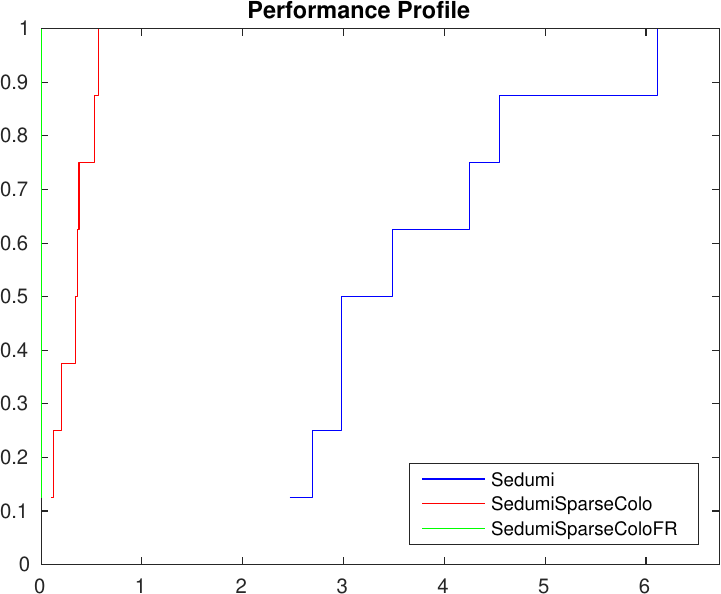}
&
\includegraphics[scale=0.3]{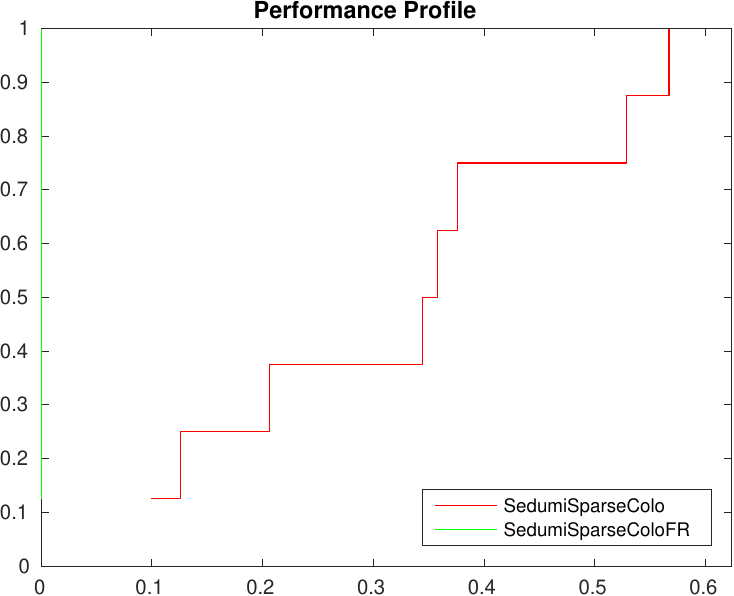}\\
Dual & 
(Dual, zoomed-in)
\end{tabular}
\end{center}
\vskip 3mm
\caption{\label{sedumivisionpp} Performance profiles (i.e., proportion of instances from a test set solved over time, relative to the best performer, cf. \cite{Dolan2002}) comparing SeDuMi with no (Sedumi), chordal-decomposition (SedumiSparseColo), and two-step processing (SedumiSparseColoFR)  on instances from polynomial optimisation.}
\end{figure}

Next, we illustrate the results on instances from a well-known polynomial optimisation problem (POP). 
In particular, the so called Lavaei-Low relaxation \cite{6756976} of the alternating-current optimal power flows (ACOPF)
has been shown  \cite{Marecek2016}  to coincide with the first level of the moment-SOS hierarchy \cite{Anjos2012} for POP. 
Due to the fact that real-life electricity transmission systems are sparse, there is sparsity present in the instances as well,
which is widely solved with SparseCoLo pre-processing or related methods \cite{Marecek2016,6692873,6756976}.
% \subsubsection{SeDuMi Results}
In our test, we consider the well-known IEEE test systems. In the name of the instance, case$x$ denotes a test system on $x$ buses, with more than $4x^2$ elements in the moment matrix, P denotes the primal SDP, and D denotes the dual SDP.
Tab.~\ref{sedumivision} presents the wall-clock run-time (including pre-processing) for 16 such SDP instances.
For larger instances (case118 and case300), the two-step pre-processing yields about 2 orders of magnitude of improvement. 
This is further illustrated by 
performance profiles on the two sets of problems
in Fig.~\ref{sedumivisionpp}:
the nearly vertical lines are for the pre-processing, while the nearly horizontal line is without the pre-processing.
This set of problems give the clearest indication of the benefits of two step pre-processing, suggesting they are particularly
structured to take advantage of the procedures.

% Note that publication-quality tables \emph{do not contain vertical rules.} 

\begin{table}
\caption{ Run-time in seconds of SeDuMi with no, chordal-completion (w/ SC), and two-step (w/ SC+FR) processing on instances from polynomial optimisation. Above mid-rule, the instances are in the primal form, below in their dual form.}
\vskip 3mm
\label{sedumivision}
\centering
\begin{tabular}{l d  d  d } 
\toprule
Instance 
 & 
\multicolumn{1}{r}{SeDuMi} & 
\multicolumn{1}{r}{\quad w/ SC \quad} & 
\multicolumn{1}{r}{w/ SC+FR} \\  
%\cmidrule(r){1-5}
\midrule
case9 &  3.47e^{+02}  &  6.40  &  4.59  \\
case14 &   2.74e^{+01}  &   2.46  &  1.63  \\
case30 &   3.36e^{+01}  &   1.49  &  1.10  \\
% 30Q is not part of the official benchmark
%case30Q & P &  9.66e+00 &  1.45e+00 & 1.29e+00 \\
case39 &   3.04e^{+01}  &   2.21  &  2.02  \\
case57 &   5.18  &   4.09e^{-01}  &  4.04e^{-01}  \\
case118 &    8.10e^{+03}  &   2.06e^{+01}  &  1.35e^{+01}  \\
case300 &   6.77e^{+01}  &   1.39  &  1.10  \\ 
\midrule
case9 &   2.57e^{+02}  &   4.83  &  2.74  \\
case14 &   2.92e^{+01}  &   2.19  &  1.99  \\
case30 &  2.63e^{+01}  &  1.16  &  8.08e^{-01}  \\
% 30Q is not part of the official benchmark
%case30Q & D &  9.09e+00 &  2.10e+00 & 1.09e+00 \\
case39 &   3.01e^{+01}  &   1.89  &  1.53  \\
case57 &   6.62  &   4.75e^{-01}  &  3.37e^{-01}  \\
case118 &   4.48e^{+03}  &   1.44e^{+01}  &  9.89  \\
case300 &   4.50e^{+01}  &   1.09  &  6.45e^{-01}  \\
\bottomrule
\end{tabular}
\end{table}

%\newcolumntype{d}{D{.}{.}{-1} } %decimal column as before

% Note that publication-quality tables \emph{do not contain vertical rules.} 

\begin{table}
\caption{\label{sedumimaxcut} Run-time in seconds of SeDuMi with no, chordal-completion (w/ SC), and two-step (w/ SC+FR) processing  on instances of the maximum cut problem from the BiqMac benchmark, with and without pre-processing.}
\vskip 3mm
\centering
\begin{tabular}{l c d  d  d } 
\toprule
N & m
 & 
\multicolumn{1}{r}{IPM} & 
\multicolumn{1}{r}{\quad w/ SC \quad} & 
\multicolumn{1}{r}{w/ SC+FR} \\  
%\cmidrule(r){1-5}
\midrule
  50 & 50 & 2.22 &  7.49e^{-01} &  3.38e^{-01} \\ 
  100 & 100 & 4.74e^{-01} & F & F \\ 
  150 & 150 &  7.87e^{-01}  & F & F \\ 
  200 & 200 &  1.30  & F & F \\ 
  300 & 300 &  1.44e^{+01}  &  8.52  &   7.56  \\ 
  350 & 350 &  1.20e^{+01}  &  1.18e^{+01}  &   1.18e^{+01}  \\ 
  400 & 400 &  1.45e^{+01}  &  1.45e^{+01}  &  1.44e^{+01}  \\ 
  450 & 450 &  1.99e^{+01}  &  1.99e^{+01}  &  1.99e^{+01}  \\ 
  500 & 500 &  3.37e^{+01}  &  3.40e^{+01}  &  3.34e^{+01}  \\ 
\midrule
  250 & 250 &  4.62  &  4.58  &  4.54  \\ 
  250 & 250 &  4.63  &  4.56  &  4.52  \\ 
  250 & 250 &  4.27  &  4.26  &  4.20  \\ 
  250 & 250 &  4.77  &  4.72  &  4.65  \\ 
  250 & 250 &  4.96  &  4.95  &  4.87  \\
\bottomrule
\end{tabular}
\end{table}

\subsection{Binary quadratic programming}

Next, let us present results on perhaps the best-known SDP relaxation, that of binary quadratic programming or, equivalently, the maximum cut problem, (MAXCUT). In Tab.~\ref{sedumimaxcut}, we see runtime values for SeDuMi by itself, with matrix completion pre-processing, and with the entire two-step pre-processing procedure. 
We indicate the size of the instance from the BiqMac benchmark as well. 
We notice that the chordal-decomposition pre-processing  sometimes results in a failure of SeDuMi,
which could be attributed to numerical failures due to degeneracy.
In that case, facial reduction does not improve upon the situation.
On small instances in the BiqMac benchmark, there is a small but consistent improvement in the overall run-time, indicating that the pre-processing does improve the efficiency.
On larger instances, the differences seem more pronounced.

%\subsection{Computer Vision}
%Next, let us consider the performance of the interior point solvers with two step pre-processing on a multi-view image registration problem problems in computer vision~\cite{schweighofer2008globally}
%and a related problem hand-eye coordination in robotics \cite{6907313}.

%\subsubsection{SDPT3 Results}

\section{Conclusions}\label{sec:conc}

Many practically-relevant instances of semidefinite programming are sparse and structured. 
Traditional general-purpose implementations of exploiting the structure \cite{Fukuda2001,Kim2012,Bergman2016}
have proven difficult to apply, due to the numerical issues they introduce \cite{Raghunathan2016}.

While one could try to exploit the structure directly, in custom code,
we suggest that general-purpose pre-processing combining both the traditional chordal-decomposition techniques 
\cite{Fukuda2001,Kim2012,Bergman2016}
and facial reduction may make it possible to exploit the structure, while relying on the
robustness of standard interior-point methods, unaware of the structure.

We have demonstrated that on SDPLib and several sets of structured problems, our combination of
chordal completion and facial reduction appears to improve the performance drastically.

\iffalse
A number of open questions remain. A thorough comparison with a fully comprehensive collection of SDP solvers could be interesting to elucidate how the effectiveness of preprocessing relates to the SDP solution method.
There is a possibility that facial reduction ameliorates the degeneracy
issues introduced by chordal-decomposition pre-processing. One could also perform a  comprehensive study numerically verifying
the presence of primal/dual non-degeneracy and constraint qualification (CQ) for the same sets of instances.
One could also ask whether similar effects are present for other classes of problems (e.g., the doubly non-negative cone and its applications in co-positive optimisation), or with solvers other than interior-point methods (e.g., first-order methods).
Finally, one could exploit recent advances in facial reduction \cite{pataki2017sieve} and 
implementation of {\sc embed} \cite{Cao2018} 
to speed up the process. 
\fi

\section{ACKNOWLEDGEMENTS}
Support for this work was provided by the OP VVV project CZ.02.1.01/0.0/0.0/16\_019/0000765 ``Research Center for Informatics"

\bibliographystyle{IEEEtran}
\bibliography{refs}

\end{document}